\documentclass[12pt]{article}%
\usepackage{amsmath}
\usepackage{amsfonts}
\usepackage{amssymb}
\usepackage{graphicx}%
\newenvironment{proof}[1][Proof]{\noindent{\textbf {#1}  }}  {\hfill$\Box$\bigskip}

\newtheorem{theorem}{Theorem}

\newtheorem{lemma}{Lemma}

\begin{document}

\title{\textbf{Large Generalized Books are $\boldsymbol{p}$-Good}}
\author{V. Nikiforov and C. C. Rousseau\\Department of Mathematical Sciences \\The University of Memphis \\373 Dunn Hall \\Memphis, Tennessee 38152-3240}
\date{}
\maketitle

\begin{abstract}
Let $B_{q}^{(r)}=K_{r}+qK_{1}$ be the graph consisting of $q$ distinct
$\left(  r+1\right)  $-cliques sharing a common $r$-clique. We prove that if
$p\geq2$ and $r\geq3$ are fixed, then
\[
r(K_{p+1},B_{q}^{(r)})=p(q+r-1)+1
\]
for all sufficiently large $q$.

\end{abstract}

Keywords: Ramsey numbers; p-good; generalized books; Szemer\'{e}di lemma

\setlength{\baselineskip}{.27in}

\section{Introduction}

The title of this paper refers to the notion of goodness introduced by Burr
and Erd\H{o}s in \cite{BE:gen} and subsequently studied by Burr and various
collaborators. A connected graph $H$ is $p$-good if the Ramsey number
$r(K_{p},H)$ is given by
\[
r(K_{p},H)=(p-1)(|V(H)|-1)+1.
\]
In this paper we prove that for every $p\geq3$ the generalized book
$B_{q}^{(r)}=K_{r}+qK_{1}$ is $p$-good if $q$ is sufficiently large.

As much as possible, standard notation is used; see, for example, \cite{Bo}. A
set of cardinality $p$ is called a $p$-set. Unless explicitly stated, all
graphs are defined on the vertex set $[n]=\{1,2,\ldots,n\}$. Let $u$ be any
vertex; then $N_{G}(u)$ and $d_{G}(u)=|N_{G}(u)|$ denote its neighborhood and
degree respectively. A graph with $n$ vertices and $m$ edges will be
designated by $G(n,m)$. By an $r$-\emph{book} we shall mean some number of
independent vertices that are each connected to every vertex of an $r$-clique.
The given $r$-clique is called the \emph{base} of the $r$-book and the
additional vertices are called the \emph{pages}. The number of pages of an
$r$-book is called its \emph{size}; the size of the largest $r$-book in a
graph $G$ is denoted by $bs^{(r)}(G)$. We shall denote the complete
$p$-partite graph with each part having $q$ vertices by $K_{p}(q)$. The
\emph{Ramsey number} $r(H_{1},H_{2})$ is the least number $n$ such that for
every graph $G$ of order $n$ either $H_{1}\subset G$ or $H_{2}\subset
\overline{G}$.

\section{The structure of subsaturated $K_{p+1}$-free graphs}

We shall need the following theorem of Andr\'{a}sfai, Erd\H{o}s and S\'{o}s
\cite{AES}.

\begin{theorem}
\label{tAES} If $G$ is a $K_{p+1}$-free graph of order $n$ and
\[
\delta(G)>\left(  1-\frac{3}{3p-1}\right)  n,
\]
then $G$ is $p$-chromatic.\hfill$\square$
\end{theorem}

The celebrated theorem of Tur\'{a}n gives a tight bound on the maximum size of
a $K_{p}$-free graph of given order. In the following theorem we show that if
the size of a $K_{p+1}$-free graph is close to the maximum then we may delete
a small portion of its vertices so that the remaining graph is $p$-chromatic.
This is a particular stability theorem in extremal graph theory (see \cite{Si}).

\begin{theorem}
\label{th1} For every $p\geq2$ there exists $c=c\left(  p\right)  >0,$ such
that for every $\alpha$ satisfying $0<\alpha\leq c,$ every $K_{p+1}$-free
graph $G=G\left(  n,m\right)  $ satisfying
\[
m\geq\left(  \frac{p-1}{2p}-\alpha\right)  n^{2}%
\]
contains an induced $p$-chromatic graph $G_{0}$ of order at least $\left(
1-2\alpha^{1/3}\right)  n$ and with minimum degree
\[
\delta\left(  G_{0}\right)  \geq\left(  1-\frac{1}{p}-4\alpha^{1/3}\right)
n.
\]

\end{theorem}

\begin{proof}
Let $c_{0}$ be the smallest positive root of the equation
\begin{equation}
x^{3}+\left(  1+\frac{3}{3p-1}\left(  \frac{p-1}{p}\right)  ^{2}\right)
x-\frac{1}{2\left(  3p-1\right)  p}=0 \label{ceq}%
\end{equation}
and set $c\left(  p\right)  =c_{0}^{3};$ then, for every $y$ satisfying
$0<y\leq c\left(  p\right)  ,$ we easily see that%
\begin{equation}
y+\left(  1+\frac{3}{3p-1}\left(  \frac{p-1}{p}\right)  ^{2}\right)
y^{1/3}\leq\frac{1}{2\left(  3p-1\right)  p}. \label{ineq01}%
\end{equation}

A rough approximation of the function $c\left(  p\right)  $ is $c(p)\approx
6^{-3}p^{-6}$, obtained by neglecting the $x^{3}$ term in equation (\ref{ceq})
and substituting the appropriate asymptotic (for large $p$) approximations for
the remaining coefficients. This gives reasonable values even for small $p$.
For all $p\geq2$,
\begin{equation}
\frac{1}{\left(  2p(3p+2)\right)  ^{3}}<c(p)<\frac{1}{\left(  2p(3p-1)\right)
^{3}}. \label{bndc}%
\end{equation}
The upper bound is evident, and the lower bound follows from a simple computation.

Let $0<\alpha\leq c\left(  p\right)  $ and the graph $G=G\left(  n,m\right)  $
satisfy the hypothesis of the theorem. We shall prove first that
\begin{equation}
\sum_{u=1}^{n}d^{2}\left(  u\right)  \leq2\left(  \frac{p-1}{p}\right)  mn.
\label{ineq1}%
\end{equation}

Indeed, writing $k_{3}(G)$ for the number of triangles in $G,$ we have
\[
3k_{3}(G)=\sum_{uv\in E}\left\vert N(u)\cap N(v)\right\vert \geq\sum_{uv\in
E}(d(u)+d(v)-n)=\sum_{u=1}^{n}d^{2}(u)-mn.
\]
Applying Tur\'{a}n's theorem to the $K_{p}$-free neighborhoods of vertices of
$G$, we deduce
\[
3k_{3}(G)\leq\frac{p-2}{2(p-1)}\sum_{u=1}^{n}d^{2}(u).
\]
Hence,
\[
\sum_{u=1}^{n}d^{2}(u)-mn\leq\frac{p-2}{2(p-1)}\sum_{u=1}^{n}d^{2}(u)
\]
and (\ref{ineq1}) follows.

Since $0<\alpha\leq c\left(  p\right)  ,$ taking the upper bound in
(\ref{bndc}) for $p=2,$ we see that $\alpha\leq20^{-3}.$ Hence,%
\begin{align*}
\left(  1+8\alpha\right)  \frac{4m^{2}}{n}  &  \geq2\left(  1+8\alpha\right)
\left(  \frac{p-1}{p}-2\alpha\right)  mn\\
&  =2\left(  \frac{p-1}{p}+\left(  6-\frac{8}{p}\right)  \alpha-16\alpha
^{2}\right)  mn\\
&  \geq2\left(  \frac{p-1}{p}+2\alpha-16\alpha^{2}\right)  mn>2\left(
\frac{p-1}{p}\right)  mn,
\end{align*}
and from (\ref{ineq1}) we deduce%
\begin{align}
\sum_{u=1}^{n}\left(  d\left(  u\right)  -\frac{2m}{n}\right)  ^{2}  &
=\sum_{u=1}^{n}d^{2}\left(  u\right)  -\frac{4m}{n}^{2}\leq2\left(  \frac
{p-1}{p}\right)  mn-\frac{4m}{n}^{2}\nonumber\\
&  <8\alpha\frac{4m^{2}}{n}\leq8\alpha\left(  \frac{p-1}{p}\right)  ^{2}n^{3}.
\label{bound3}%
\end{align}

Set $V=V\left(  G\right)  $ and let $M_{\varepsilon}$ be the set of all
vertices $u\in V$ satisfying $d\left(  u\right)  <2m/n-\varepsilon n.$ For
every $\varepsilon>0,$ inequality (\ref{bound3}) implies
\[
\left\vert M_{\varepsilon}\right\vert \varepsilon^{2}n^{2}<\sum_{u\in
M_{\varepsilon}}\left(  d\left(  u\right)  -\frac{2m}{n}\right)  ^{2}\leq
\sum_{u\in V}\left(  d\left(  u\right)  -\frac{2m}{n}\right)  ^{2}\leq
8\alpha\left(  \frac{p-1}{p}\right)  ^{2}n^{3},
\]
and thus,
\begin{equation}
\left\vert M_{\varepsilon}\right\vert <8\varepsilon^{-2}\alpha\left(
\frac{p-1}{p}\right)  ^{2}n. \label{bound4.0}%
\end{equation}

Furthermore, setting $G_{\varepsilon}=G\left[  V\backslash M_{\varepsilon
}\right]  ,$ for every $u\in V\left(  G_{\varepsilon}\right)  ,$ we obtain%
\begin{equation}
d_{G_{\varepsilon}}\left(  u\right)  \geq d\left(  u\right)  -\left\vert
M_{\varepsilon}\right\vert \geq\frac{2m}{n}-\varepsilon n-\left\vert
M_{\varepsilon}\right\vert >\frac{p-1}{p}n-2\alpha n-\varepsilon n-\left\vert
M_{\varepsilon}\right\vert . \label{bound4.1}%
\end{equation}

For $\varepsilon=2\alpha^{1/3}$ we claim that%
\begin{equation}
\frac{p-1}{p}n-2\alpha n-\varepsilon n-\left\vert M_{\varepsilon}\right\vert
>\frac{3p-4}{3p-1}\left(  n-\left\vert M_{\varepsilon}\right\vert \right)
=\frac{3p-4}{3p-1}v\left(  G_{\varepsilon}\right)  . \label{bound4}%
\end{equation}
Indeed, assuming the opposite and applying inequality (\ref{bound4.0}) with
$\varepsilon=2\alpha^{1/3},$ we see that%
\[
\left(  \frac{1}{\left(  3p-1\right)  p}-2\alpha-2\alpha^{1/3}\right)
n\leq\frac{3}{3p-1}\left\vert M_{2\alpha^{1/3}}\right\vert <2\frac{3}%
{3p-1}\left(  \frac{p-1}{p}\right)  ^{2}\alpha^{1/3}n;
\]
hence,
\[
2\alpha+2\left(  1+\frac{3}{3p-1}\left(  \frac{p-1}{p}\right)  ^{2}\right)
\alpha^{1/3}-\frac{1}{\left(  3p-1\right)  p}>0,
\]
contradicting (\ref{ineq01}).

Set $G_{0}=G_{2\alpha^{1/3}}$; from (\ref{bound4}), we see that $G_{0}$
satisfies the conditions of Theorem \ref{tAES}, so it is $p$-chromatic.

Finally, from (\ref{bound4.0}) and (\ref{bound4.1}), we have
\begin{align*}
\delta\left(  G_{0}\right)   &  \geq\frac{p-1}{p}n-2\alpha n-2\alpha
^{1/3}n-\left(  \frac{p-1}{p}\right)  ^{2}\alpha^{1/3}n>\frac{p-1}{p}n-2\alpha
n-3\alpha^{1/3}n\\
&  >\left(  1-\frac{1}{p}-4\alpha^{1/3}\right)  n,
\end{align*}
completing the proof.
\end{proof}

\section{A Ramsey property of ${\boldsymbol{K}_{p+1}}$-free graphs}

The main result of this section is the following theorem.

\begin{theorem}
\label{th2} Let $r\geq2,$ $p\geq2$ be fixed. For every $\xi>0$ there exists an
$n_{0}=n_{0}(p,r,\xi)$ such that every graph $G$ of order $n\geq n_{0}$ that
is $K_{p+1}$-free either satisfies $bs^{(r)}\left(  \overline{G}\right)
>n/p$, or contains an induced $p$-chromatic graph $G_{1}$ of order $\left(
1-\xi\right)  n$ and minimum degree
\[
\delta(G_{1})\geq\left(  1-\frac{1}{p}-2\xi\right)  n.
\]

\end{theorem}

Our main tool in the proof of Theorem \ref{th2} is the regularity lemma of
Szemer\'{e}di (SRL for short); for expository matter on SRL see \cite{Bo} and
\cite{KoSi}. For the sake of completeness we formulate here the relevant basic notions.

Let $G$ be a graph; if $A,B\subset V\left(  G\right)  $ are nonempty disjoint
sets, we write $e\left(  A,B\right)  $ for the number of $A-B$ edges and call
the value
\[
d\left(  A,B\right)  =\frac{e\left(  A,B\right)  }{\left\vert A\right\vert
\left\vert B\right\vert }%
\]
the \emph{density} of the pair $\left(  A,B\right)  .$

Let $\varepsilon>0;$ a pair $\left(  A,B\right)  $ of two nonempty disjoint
sets $A,B\subset V\left(  G\right)  $ is called $\varepsilon$\emph{-regular}
if the inequality%
\[
\left\vert d\left(  A,B\right)  -d\left(  X,Y\right)  \right\vert <\varepsilon
\]
holds whenever $X\subset A,$ $Y\subset B,$ $\left\vert X\right\vert
\geq\varepsilon\left\vert A\right\vert ,$ and $\left\vert Y\right\vert
\geq\varepsilon\left\vert B\right\vert .$

We shall use SRL in the following form.

\begin{theorem}
[Szemer\'{e}di's Regularity Lemma]\label{SRL} Let $l\geq1$, $\varepsilon>0$.
There exists $M=M\left(  \varepsilon,l\right)  $ such that, for every graph
$G$ of sufficiently large order $n$, there exists a partition $V\left(
G\right)  =\cup_{i=0}^{k}V_{i}$ satisfying $l\leq k\leq M$ and:

$\emph{(i)}$ $\left\vert V_{0}\right\vert <\varepsilon n,$ $\left\vert
V_{1}\right\vert =...=\left\vert V_{k}\right\vert ;$

\emph{(ii)} all but at most $\varepsilon k^{2}$ pairs $\left(  V_{i}%
,V_{j}\right)  ,$ $\left(  i,j\in\left[  k\right]  \right)  ,$ are
$\varepsilon$-uniform.
\end{theorem}

We also need a few technical results; the first one is a basic property of
$\varepsilon$-regular pairs (see \cite{KoSi}, Fact 1.4).

\begin{lemma}
\label{XPle} Suppose $0<\varepsilon<d\leq1$ and $\left(  A,B\right)  $ is an
$\varepsilon$-regular pair with $e\left(  A,B\right)  =d|A||B|$. If $Y\subset
B$ and $(d-\varepsilon)^{r-1}|Y|>\varepsilon|B|$ where $r>1$, then there are
at most $\varepsilon r|A|^{r}$ $r$-sets $R\subset A$ with
\[
\left\vert \left(  \bigcap_{u\in R}N\left(  u\right)  \right)  \cap
Y\right\vert \leq(d-\varepsilon)^{r}|Y|.
\]

\end{lemma}

The next lemma gives a lower bound on the number of $r$-cliques in a graph
consisting of several dense $\varepsilon$-regular pairs sharing a common part.

\begin{lemma}
\label{dle} Suppose $0<\varepsilon<d\leq1$ and $(d-\varepsilon)^{r-2}%
>\varepsilon$. Suppose $H$ is a graph and $V(H)=A\cup B_{1}\cup\cdots\cup
B_{t}$ is a partition with $\left\vert A\right\vert =|B_{1}|=\cdots=|B_{t}|$
and such that for every $i\in\lbrack t]$ the pair $(A,B_{i})$ is $\varepsilon
$-regular with $e(A,B_{i})\geq d|A||B_{i}|$. If $m$ is the number of the
$r$-cliques in $A,$ then at least
\[
t|A|\left(  m-\varepsilon r|A|^{r}\right)  (d-\varepsilon)^{r}%
\]
$(r+1)$-cliques of $H$ have exactly $r$ vertices in $A$.
\end{lemma}

\begin{proof}
Set $a=|A|=|B_{1}|=\cdots=|B_{t}|$. For every $i\in\left[  t\right]  ,$
applying Lemma \ref{XPle} to the pair $(A,B_{i})$ with $Y=B_{i}$ we conclude
that there are at most $\varepsilon ra^{r-1}$ $r$-sets $R\subset A$ with
\[
\left\vert \left(  \bigcap_{u\in R}N\left(  u\right)  \right)  \cap
B_{i}\right\vert \leq\left(  d-\varepsilon\right)  ^{r}a,
\]
and therefore, at least $(m-\varepsilon ra^{r})$ $r$-cliques $R\subset A$
satisfy
\[
\left\vert \left(  \bigcap_{u\in R}N\left(  u\right)  \right)  \cap
B_{i}\right\vert >(d-\varepsilon)^{r}a.
\]
Hence, at least $t(d-\varepsilon)^{r}(m-\varepsilon ra^{r})a$ $(r+1)$-cliques
of $H$ have exactly $r$ vertices in $A$ and one vertex in $\cup_{i\in\left[
t\right]  }B_{i},$ completing the proof.
\end{proof}

The following consequence of Ramsey's theorem has been proved by Erd\H{o}s
\cite{Erd}.

\begin{lemma}
\label{le2} Given integers $p \ge2$, $r \ge2$, there exist a $c_{p,r}>0$ such
that if $G$ is a $K_{p+1}$-free graph of order $n$ and $n\geq r(K_{p+1}%
,K_{r})$ then $G$ contains at least $c_{p,r}n^{r}$ independent $r$-sets.
\end{lemma}

We need another result related to the regularity lemma of Szemer\'{e}di, the
so-called Key Lemma (e.g., see \cite{KoSi}, Theorem 2.1). We shall use the
following simplified version of the Key Lemma.

\begin{theorem}
\label{tKL} Suppose $0<\varepsilon<d<1$ and let $m$ be a positive integer. Let
$G$ be a graph of order $(p+1)m$ and let $V(G)=V_{1}\cup\cdots\cup V_{p+1}$ be
a partition of $V(G)$ into $p+1$ sets of cardinality $m$ so that each of the
pairs $(V_{i},V_{j})$ is $\varepsilon$-regular and has density at least $d$.
If $\varepsilon\leq(d-\varepsilon)^{p}/(p+2)$ then $K_{p+1}\subset G$.
\end{theorem}

\begin{proof}
[Proof of Theorem \ref{th2}]\textbf{\ }Our proof is straightforward but rather
rich in technical details, so we shall briefly outline it first. For some
properly selected $\varepsilon$, applying SRL, we partition all but
$\varepsilon n$ vertices of $G$ in $k$ sets $V_{1},\ldots,V_{k}$ of equal
cardinality such that almost all pairs $(V_{i},V_{j})$ are $\varepsilon
$-regular. We may assume that the number of dense $\varepsilon$-regular pairs
$(V_{i},V_{j})$ is no more than $\frac{p-1}{2p}k^{2},$ since otherwise, from
Theorem \ref{tKL} and Tur\'{a}n's theorem, $G$ will contain a $K_{p+1}$.
Therefore, there are at least $(1/2p+o(1))k^{2}$ sparse $\varepsilon$-regular
pairs $(V_{i},V_{j})$. From Lemma \ref{le2} it follows that the number of
independent $r$-sets in any of the sets $V_{1},\ldots,V_{k}$ is $\Theta
(n^{r})$. Consider the size of the $r$-book in $\overline{G}$ having for its
base the average independent $r$-set in $V_{i}$. For every sparse
$\varepsilon$-regular pair $(V_{i},V_{j})$ almost every vertex in $V_{j}$ is a
page of such a book. Also each $\varepsilon$-regular pair $(V_{i},V_{j})$
whose density is not very close to 1 contributes substantially many additional
pages to such books. Precise estimates show that either $bs^{(r)}(\overline
{G})>n/p$ or else the number of all $\varepsilon$-regular pairs $(V_{i}%
,V_{j})$ with density close to $1$ is $\left(  \frac{p-1}{2p}+o(1)\right)
k^{2}$. Thus the size of $G$ is $\left(  \frac{p-1}{2p}+o(1)\right)  n^{2}$
and therefore, according to Theorem \ref{th1}, $G$ contains the required
induced $p$-chromatic subgraph with the required minimum degree.

\textit{Details of the proof.} Let $c(p)$ be as in Theorem \ref{th1} and
$c_{p,r}$ be as in Lemma \ref{le2}. Select
\begin{equation}
\delta=\min\left\{  \frac{\xi^{3}}{32},\frac{c(p)}{4}\right\}  , \label{gdd}%
\end{equation}
set
\begin{equation}
d=\min\left\{  \left(  \frac{\delta}{2}\right)  ^{r+1}\left(  \frac{r}%
{c_{p,r}}+2r+1+2p\right)  ^{-1},\;\frac{p\delta}{1+p\delta}\left(  \frac
{r}{c_{p,r}}+2r+1\right)  ^{-1}\right\}  , \label{ddef}%
\end{equation}
and let
\begin{equation}
\varepsilon=\min\left\{  \delta,\frac{d^{p}}{2\left(  p+1\right)  }\right\}  .
\label{edef}%
\end{equation}

These definitions are justified at the later stages of the proof. Since
$c_{p,r}<r!$ we easily see that $0<2\varepsilon<d<\delta<1.$ Hence,
Bernoulli's inequality implies
\begin{equation}
(d-\varepsilon)^{p}\geq d^{p}-p\varepsilon d^{p-1}>d^{p}-p\varepsilon
=2(p+1)\varepsilon-p\varepsilon=(p+2)\varepsilon. \label{dbnd}%
\end{equation}

Applying SRL we find a partition $V\left(  G\right)  =V_{0}\cup V_{1}%
\cup\cdots\cup V_{k}$ so that $\left\vert V_{0}\right\vert <\varepsilon n$,
$\left\vert V_{1}\right\vert =\cdots=\left\vert V_{k}\right\vert $ and all but
$\varepsilon k^{2}$ pairs $(V_{i},V_{j})$ are $\varepsilon$-regular. Without
loss of generality we may assume $\left\vert V_{i}\right\vert >r(K_{p+1}%
,K_{r})$ and $k>1/\varepsilon$. Consider the graphs $H_{\text{irr}},$
$H_{\text{lo}}$, $H_{\text{mid}}$ and $H_{\text{hi}}$ defined on the vertex
set $[k]$ as follows:

\begin{itemize}
\item[(i)] $(i,j) \in E(H_{\text{irr}}) $ iff the pair $(V_{i}, V_{j})$ is not
$\varepsilon$-regular,

\item[(ii)] $(i,j)\in E(H_{\text{lo}})$ iff the pair $(V_{i},V_{j})$ is
$\varepsilon$-regular and
\[
d(V_{i},V_{j})\leq d,
\]

\item[(iii)] $(i,j)\in E(H_{\text{mid}})$ iff the pair $(V_{i},V_{j})$ is
$\varepsilon$-regular and
\[
d<d(V_{i},V_{j})\leq1-\delta,
\]

\item[(iv)] $(i,j)\in E(H_{\text{hi}})$ iff the pair $(V_{i},V_{j})$ is
$\varepsilon$-regular and
\[
d(V_{i},V_{j})>1-\delta.
\]

\end{itemize}

Clearly, no two of these graphs have edges in common; thus
\[
e(H_{\text{irr}})+e(H_{\text{lo}})+e(H_{\text{mid}})+e(H_{\text{hi}}%
)=\binom{k}{2}.
\]
Hence, from $d>2\varepsilon$ and $k>1/\varepsilon,$ we see that
\begin{align}
e(H_{\text{lo}})+e(H_{\text{mid}})+e(H_{\text{hi}})  &  \geq\binom{k}%
{2}-\varepsilon k^{2}=\frac{k^{2}}{2}-\frac{k}{2}-\varepsilon k^{2}\nonumber\\
&  \geq\frac{k^{2}}{2}-\varepsilon k^{2}-\varepsilon k^{2}>\left(  \frac{1}%
{2}-d\right)  k^{2}. \label{ineq3}%
\end{align}

Since $G$ is $K_{p+1}$-free, from (\ref{dbnd}), we have $\varepsilon
\leq\left(  d-\varepsilon\right)  ^{p}/(p+2);$ applying Theorem \ref{tKL}, we
conclude that the graph $H_{\text{mid}}\cup H_{\text{hi}}$ is $K_{p+1}$-free.
Therefore, from Tur\'{a}n's theorem,
\[
e(H_{\text{mid}})+e(H_{\text{hi}})\leq\left(  \frac{p-1}{2p}\right)  k^{2},
\]
and from inequality (\ref{ineq3}) we deduce
\begin{equation}
e(H_{\text{lo}})>\left(  \frac{1}{2p}-d\right)  k^{2}. \label{ineq2}%
\end{equation}

Next we shall bound $bs^{(r)}(\overline{G})$ from below. To achieve this we
shall count the independent $(r+1)$-sets having exactly $r$ vertices in some
$V_{i}$ and one vertex outside $V_{i}$. Fix $i\in\lbrack k]$ and let $m$ be
the number of independent $r$-sets in $V_{i}$. Observe that Lemma \ref{le2}
implies $m\geq c_{p,r}|V_{i}|^{r}$.

Set $L=N_{H_{\text{lo}}}(i)$ and apply Lemma \ref{dle} with $A=V_{i}$,
$B_{j}=V_{j},$ for all $j\in L$, and
\[
H=\overline{G}\left[  A\cup\left(  \bigcup_{j\in L}B_{j}\right)  \right]  .
\]
Since, for every $j\in L,$ the pair $(V_{i},V_{j})$ is $\varepsilon$-regular
and
\[
e_{H}(V_{i},V_{j})\geq(1-d)|V_{i}||V_{j}|,
\]
we conclude that there are at least%
\[
d_{H_{\text{lo}}}(i)|V_{i}|(m-\varepsilon r|V_{i}|^{r})(1-d-\varepsilon)^{r}%
\]
independent $(r+1)$-sets in $G$ having exactly $r$ vertices in $V_{i}$ and one
vertex in $\cup_{j\in L}B_{j}$.

Set now $M=N_{H_{\text{mid}}}(i)$, and apply Lemma \ref{dle} with $A=V_{i}$,
$B_{j}=V_{j}$ for all $j\in M$ and
\[
H=\overline{G}\left[  A\cup\left(  \bigcup_{j\in M}B_{j}\right)  \right]  .
\]
Since, for every $j\in M,$ the pair $(V_{i},V_{j})$ is $\varepsilon$-regular
and
\[
e_{H}(V_{i},V_{j})\geq\delta|V_{i}||V_{j}|,
\]
we conclude that there are at least
\[
d_{H_{\text{mid}}}(i)|V_{i}|\left(  m-\varepsilon r|V_{i}|^{r}\right)
(\delta-\varepsilon)^{r}%
\]
independent $(r+1)$-sets in $G$ having exactly $r$ vertices in $V_{i}$ and one
vertex in $\cup_{j\in L}B_{j}$. Since
\[
\left(  \bigcup_{j\in L}B_{j}\right)  \bigcap\left(  \bigcup_{j\in M}%
B_{j}\right)  =\varnothing,
\]
there are at least
\[
d_{H_{\text{lo}}}(i)|V_{i}|\left(  m-\varepsilon r|V_{i}|^{r}\right)
(1-d-\varepsilon)^{r}+d_{H_{\text{mid}}}(i)|V_{i}|\left(  m-\varepsilon
r|V_{i}|^{r}\right)  (\delta-\varepsilon)^{r}%
\]
independent $(r+1)$-sets in $G$ having exactly $r$ vertices in $V_{i}$ and one
vertex outside $V_{i}$. Thus, taking the average over all $m$ independent
$r$-sets in $V_{i}$, we conclude
\begin{align*}
bs^{\left(  r\right)  }\left(  \overline{G}\right)   &  \geq\left\vert
V_{i}\right\vert \left(  1-\frac{\varepsilon r}{c_{p,r}}\right)  \left(
d_{H_{\text{lo}}}\left(  i\right)  \left(  1-d-\varepsilon\right)
^{r}+d_{H_{\text{mid}}}\left(  i\right)  \left(  \delta-\varepsilon\right)
^{r}\right) \\
&  \geq n\left(  \frac{1-\varepsilon}{k}\right)  \left(  1-\frac{\varepsilon
r}{c_{p,r}}\right)  \left(  d_{H_{\text{lo}}}\left(  i\right)  \left(
1-d-\varepsilon\right)  ^{r}+d_{H_{\text{mid}}}\left(  i\right)  \left(
\delta-\varepsilon\right)  ^{r}\right)  .
\end{align*}
Summing this inequality for all $i=1,\ldots,k$ we obtain%
\begin{align}
\frac{bs^{\left(  r\right)  }\left(  \overline{G}\right)  }{n}  &  \geq\left(
1-\varepsilon\right)  \left(  1-\frac{\varepsilon r}{c_{p,r}}\right)  \left(
\frac{2e\left(  H_{\text{lo}}\right)  }{k^{2}}\left(  1-d-\varepsilon\right)
^{r}+\frac{2e\left(  H_{\text{mid}}\right)  }{k^{2}}\left(  \delta
-\varepsilon\right)  ^{r}\right) \nonumber\\
&  >\left(  1-\left(  \frac{r}{c_{p,r}}+1\right)  \varepsilon\right)  \left(
\frac{2e\left(  H_{\text{lo}}\right)  }{k^{2}}\left(  1-r\left(
d+\varepsilon\right)  \right)  +\frac{2e\left(  H_{\text{mid}}\right)  }%
{k^{2}}\left(  \delta-\varepsilon\right)  ^{r}\right) \nonumber\\
&  >\left(  1-\left(  \frac{r}{c_{p,r}}+1\right)  d\right)  \left(
\frac{2e\left(  H_{\text{lo}}\right)  }{k^{2}}\left(  1-2rd\right)
+\frac{2e\left(  H_{\text{mid}}\right)  }{k^{2}}\left(  \frac{\delta}%
{2}\right)  ^{r}\right) \nonumber\\
&  >\left(  1-\left(  \frac{r}{c_{p,r}}+2r+1\right)  d\right)  \frac{2e\left(
H_{\text{lo}}\right)  }{k^{2}}+\left(  1-\left(  \frac{r}{c_{p,r}}+1\right)
d\right)  \left(  \frac{\delta}{2}\right)  ^{r}\frac{2e\left(  H_{\text{mid}%
}\right)  }{k^{2}}. \label{in4}%
\end{align}

Assume the assertion of the theorem false and suppose
\begin{equation}
bs^{\left(  r\right)  }\left(  \overline{G}\right)  \leq\frac{n}{p}.
\label{bsb}%
\end{equation}

We shall prove that this assumption implies
\begin{align}
e\left(  H_{\text{lo}}\right)   &  <\left(  \frac{1}{2p}+\frac{\delta}%
{2}\right)  k^{2},\label{in1}\\
e\left(  H_{\text{mid}}\right)   &  <\delta k^{2}. \label{in2}%
\end{align}

Disregarding the term $e\left(  H_{\text{mid}}\right)  $ in (\ref{in4}), in
view of (\ref{bsb}) and (\ref{ddef}), we have
\begin{align*}
e\left(  H_{\text{lo}}\right)   &  <\left(  1-\left(  \frac{r}{c_{p,r}%
}+2r+1\right)  d\right)  ^{-1}\frac{bs^{\left(  r\right)  }\left(
\overline{G}\right)  }{2n}k^{2}\\
&  \leq\left(  1-\left(  \frac{r}{c_{p,r}}+2r+1\right)  d\right)  ^{-1}%
\frac{k^{2}}{2p}\\
&  \leq\left(  1-\frac{p\delta}{1+p\delta}\right)  ^{-1}\frac{k^{2}}%
{2p}=\left(  \frac{1}{2p}+\frac{\delta}{2}\right)  k^{2},
\end{align*}
and inequality (\ref{in1}) is proved.

Furthermore, observe that equality (\ref{ddef}) implies%
\[
\left(  \frac{r}{c_{p,r}}+1\right)  d<\left(  \frac{r}{c_{p,r}}+2r+1\right)
d\leq\frac{p\delta}{1+p\delta}\leq p\delta<\frac{1}{2},
\]
and consequently,
\[
\left(  1-\left(  \frac{r}{c_{p,r}}+1\right)  d\right)  >\frac{1}{2}.
\]

Hence, from (\ref{in4}), taking into account (\ref{bsb}) and (\ref{ineq2}), we
find that
\begin{align*}
\frac{e\left(  H_{\text{mid}}\right)  }{2}\left(  \frac{\delta}{2}\right)
^{r}  &  <e\left(  H_{\text{mid}}\right)  \left(  \frac{\delta}{2}\right)
^{r}\left(  1-\left(  \frac{r}{c_{p,r}}+1\right)  d\right) \\
&  \leq\frac{bs^{\left(  r\right)  }\left(  \overline{G}\right)  k^{2}}%
{2n}-\left(  1-\left(  \frac{r}{c_{p,r}}+2r+1\right)  d\right)  e\left(
H_{\text{lo}}\right) \\
&  <\left(  \frac{1}{2p}-\left(  1-\left(  \frac{r}{c_{p,r}}+2r+1\right)
d\right)  \left(  \frac{1}{2p}-d\right)  \right)  k^{2}\\
&  =\left(  1+\left(  \frac{r}{c_{p,r}}+2r+1\right)  \left(  \frac{1}%
{2p}-d\right)  \right)  dk^{2}\\
&  <\frac{1}{2p}\left(  \frac{r}{c_{p,r}}+2r+1+2p\right)  dk^{2}<\left(
\frac{\delta}{2}\right)  ^{r+1}k^{2}.
\end{align*}
Therefore, inequality (\ref{in2}) holds also.

Furtermore, inequality (\ref{ineq3}), together with (\ref{in1}) and
(\ref{in2}), implies
\[
e\left(  H_{\text{hi}}\right)  >\left(  \frac{1}{2}-d\right)  k^{2}-\left(
\frac{1}{2p}+\frac{\delta}{2}\right)  k^{2}-\delta k^{2}=\left(  \frac
{p-1}{2p}-\frac{5\delta}{2}\right)  k^{2},
\]
and consequently, from the definition of $H_{\text{hi}},$ we obtain
\begin{align*}
e\left(  G\right)   &  \geq e\left(  H_{\text{hi}}\right)  \left(
\frac{\left(  1-\varepsilon\right)  n}{k}\right)  ^{2}\left(  1-\delta\right)
>\left(  \frac{p-1}{2p}-\frac{5\delta}{2}\right)  \left(  1-2\varepsilon
\right)  \left(  1-\delta\right)  n^{2}\\
&  =\frac{p-1}{2p}\left(  1-\frac{5p\delta}{p-1}\right)  \left(
1-2\varepsilon\right)  \left(  1-\delta\right)  n^{2}>\\
&  >\frac{p-1}{2p}\left(  1-\left(  \frac{5p}{p-1}+3\right)  \delta\right)
n^{2}>\left(  \frac{p-1}{2p}-4\delta\right)  n^{2}.
\end{align*}
Hence, by (\ref{gdd}), applying Theorem \ref{th1}, it follows that $G$
contains an induced $p$-chromatic graph with the required properties.
\end{proof}

Following the basic idea of the proof of Theorem \ref{th2} but applying the
complete Key Lemma instead of Theorem \ref{tKL}, we obtain a more general
result, whose proof, however, is considerably easier than the proof of Theorem
\ref{th2}.

\begin{theorem}
\label{th3} Suppose $H$ is a fixed $\left(  p+1\right)  $-chromatic graph. For
every $H$-free graph $G$ of order $n,$%
\[
bs^{\left(  r\right)  }\left(  \overline{G}\right)  >\left(  \frac{1}%
{p}+o\left(  1\right)  \right)  n.
\]
\hfill\ $\ \hfill\square$
\end{theorem}

Note that the graph $K_{p}\left(  q+r-1\right)  $ is $p$-chromatic and its
complement has no $B_{q}^{(r)},$ so for every $\left(  p+1\right)  $-chromatic
graph $H$ and every $r,q$ we have%
\[
r\left(  H,B_{q}^{\left(  r\right)  }\right)  \geq p\left(  q+r-1\right)  +1.
\]
Hence, from Theorem \ref{th3}, we immediately obtain the following theorem.

\begin{theorem}
\label{th4} For every fixed $\left(  p+1\right)  $-chromatic graph $H$ and
fixed integer $r>1,$%
\[
r\left(  H,B_{q}^{\left(  r\right)  }\right)  =pq+o\left(  q\right)  .
\]
\hfill\ $\ \hfill\square$
\end{theorem}

Note that it is not possible to avoid the $o\left(  q\right)  $ term in
Theorem \ref{th4} without additional stipulations about $H$, since, as
Faudree, Rousseau and Sheehan have shown in \cite{FRS}, the inequality%
\[
r\left(  C_{4},B_{q}^{(2)}\right)  \geq q+2\sqrt{q}%
\]
holds for infinitely many values of $q.$ However, when $H=K_{p+1}$ and $q$ is
large we can prove a precise result.

\section{Ramsey numbers $r\left(  K_{p},B_{q}^{\left(  r\right)  }\right)  $
for large $q$}

In this section we determine $r\left(  K_{p},B_{q}^{\left(  r\right)
}\right)  $ for fixed $p\geq3,$ $r\geq2$ and large $q.$

\begin{theorem}
For fixed $p\geq2$ and $r\geq2$, $r(K_{p+1},B_{q}^{(r)})=p(q+r-1)+1$ for all
sufficiently large $q$.
\end{theorem}

\begin{proof}
Since $K_{p}(q+r-1)$ contains no $K_{p+1}$ and its complement contains no
$B_{q}^{(r)}$, we have
\[
r(K_{p+1},B_{q}^{(r)})\geq p(q+r-1)+1.
\]

Let $G$ be a $K_{p+1}$-free graph of order $n=p(q+r-1)+1$. Since $n/p>q$,
either we're done or else $G$ contains an induced $p$-chromatic subgraph
$G_{1}$ of order $pq+o(q)$ with minimum degree
\[
\delta(G_{1})\geq\left(  1-\frac{1}{p}+o(1)\right)  n.
\]
Using this bound on $\delta(G_{1})$ we can easily prove by induction on $p$
that $G_{1}$ contains a copy of $K_{p}(r)$. Fix a copy of $K_{p}(r)$ in
$G_{1}$ and let $A_{1},A_{2},\ldots,A_{p}$ be its vertex classes. Let
$A=A_{1}\cup\cdots\cup A_{p}$ and $B=V(G)\setminus A$. If some vertex $i\in B$
is adjacent to at least one vertex in each of the parts $A_{1},A_{2}%
,\ldots,A_{p}$ then $G$ contains a $K_{p+1}$. Otherwise for each vertex $u\in
B$ there is at least one $v$ so that $u$ is adjacent in $\overline{G}$ to all
members of $A_{v}$. It follows by the pigeonhole principle that $bs^{(r)}%
(\overline{G})=s$ where
\[
s\geq\left\lceil \frac{n-p(r-1)}{p}\right\rceil =\left\lceil q-1+\frac{1}%
{p}\right\rceil =q,
\]
and we really are done.
\end{proof}

The proof using the regularity lemma that $r(K_{p+1},B_{q}^{(r)})=p\left(
q+r-1\right)  +1$ if $q$ is sufficiently large does indeed require that $q$
increase quite rapidly as a function of the parameters $p$ and $r$. This
raises the question of what growth rate is actually required. The following
simple calculation shows that polynomial growth in $p$ is not sufficient.

\begin{theorem}
For arbitrary fixed $k$ and $r$,
\[
\frac{r(K_{m},B_{m^{k}}^{(r)})}{m^{k+r-1}}\rightarrow\infty
\]
as $m\rightarrow\infty$.
\end{theorem}

\begin{proof}
We shall prove that $r(K_{m},B_{m^{k}}^{(r)})>cm^{k+r}/(\log m)^{r}$ for all
sufficiently large $m$. Let $N=\lfloor cm^{k+r}/(\log m)^{r}\rfloor$ where $c$
is to be chosen, and set $p=(C/m)\log m$ where $C=2(k+r-1)$. Let $G$ be the
random graph $G=G(N,1-p)$. The probability that $K_{m}\subset G$
\begin{align*}
{\mathbb{P}}(K_{m}\subset G)  &  \leq\binom{N}{m}(1-p)^{\binom{m}{2}}%
\leq\binom{N}{m}e^{-pm(m-1)/2}<\left(  \frac{Ne}{m}\right)  ^{m}%
e^{pm/2}m^{-\left(  k+r-1\right)  m}\\[0.1in]
&  =\left(  \frac{Ne^{1+p/2}m^{-(k+r-1)}}{m}\right)  ^{m}=o(1),\text{
\ \ }m\rightarrow\infty.
\end{align*}

To bound the probability that $B_{m^{k}}^{(r)}\subset\overline{G}$, we use the
following simple consequence of Chernoff's inequality \cite{Be87}: if
$X=X_{1}+X_{2}+\cdots+X_{n}$ where independently each $X_{i}=1$ with
probability ${\mathfrak{p}}$ and $X_{i}=0$ with probability $1-{\mathfrak{p}}$
then
\[
{\mathbb{P}}(X\geq M)\leq\left(  \frac{n{\mathfrak{p}}e}{M}\right)  ^{M}%
\]
for any $M\geq n{\mathfrak{p}}.$ Thus we find
\[
{\mathbb{P}}(B_{m^{k}}^{(r)}\subset\overline{G})\leq\binom{N}{r}%
p^{r(r-1)/2}\left(  \frac{(N-r)p^{r}e}{m^{k}}\right)  ^{m^{k}}.
\]
Since the product of the first two factors has polynomial growth in $m$, to
have ${\mathbb{P}}(B_{m^{k}}^{(r)})=o(1)$ when $m\rightarrow\infty,$ it
suffices to take $c=1/(3C^{r}),$ so that
\[
\frac{(N-r)p^{r}e}{m^{k}}\leq\frac{(c\,m^{k+r}/(\log m)^{r})((C/m)\log
m)^{r}e}{m^{k}}=\frac{e}{3},
\]
making the last factor approach 0 exponentially.
\end{proof}

\end{document}